\documentstyle{article}

\def\g{\mbox{\bf g\,}}
\def\R{\mbox{$\cal R$\,}}
\def\F{\mbox{$\cal F$}}
\def\1{{\bf 1}}
\def\ot{\otimes}
\def\id{\mbox{id}}
\def\A{\mbox{$\cal A$}}
\def\Ah{\mbox{${\cal A}_h$}}
\newcommand{\tr}{\triangleright_h}
\newcommand{\trc}{\triangleright}
\newtheorem{prop}{Proposition}

\begin{document}
\begin{titlepage}
\begin{center}
October 1997  \hfill       LMU-TPW 97-27\\
\vskip.6in

{\Large \bf q-Deforming Maps for Lie Group Covariant Heisenberg
Algebras}\footnote{Talk presented at the Fifth Wigner Symposium,
25-29 August 1997, Vienna, Germany. Submitted for the proceedings
of the Conference.}

\vskip.4in

Gaetano Fiore

\vskip.25in

{\em Sektion Physik der Ludwig-Maximilians-Universit\"at
M\"unchen\\
Theoretische Physik --- Lehrstuhl Professor Wess\\
Theresienstra\ss e 37, 80333 M\"unchen\\
Federal Republic of Germany}

{\footnotesize 
{\it e-mail: }Gaetano.Fiore \ @ \ physik.uni-muenchen.de}
\end{center}

\vskip1in
 \begin{abstract}
We briefly summarize our systematic construction
procedure of $q$-deforming maps for Lie group covariant
Weyl or Clifford algebras.
 \end{abstract}
 \end{titlepage}

\section{Introduction}

Any deformation of a Weyl or Clifford algebra can  be realized 
through a change of generators in the undeformed algebra 
\cite{mathias,ducloux}. ``$q$-Deformations''
of Weyl or Clifford algebras that were covariant under the action of 
a simple Lie algebra $\g$ are characterized by their being 
covariant under the action of the quantum group $U_h\g$, where $q=e^h$.
Here we briefly summarize our systematic construction procedure 
\cite{fiojmp,fiocmp} of
all the possible corresponding changes of generators, together with
the corresponding realizations of the $U_h\g$-action.

This paves the way \cite{fiojmp} for a physical 
interpretation of deformed
generators as ``composite operators'', functions of the
undeformed ones. For instance, if the latter act as 
creators and annihilators on
a bosonic or fermionic Fock space, then the former would act as
creators and annihilators of
some sort of ``dressed states'' in the same space.
Since there exists  \cite{fiocmp} a
basis of $\g$-invariants that depend on the undeformed generators in a
non-polynomial way, but on the deformed ones in a polynomial way,
these changes of generators might be employed
to simplify the dynamics of some $\g$-covariant quantum physical systems
based on some complicated $\g$-invariant Hamiltonian.

Let us list the essential ingredients of our construction procedure:
\begin{enumerate}

\item \g, a simple Lie algebra.
\item The cocommutative Hopf algebra $H\equiv(U\g,\cdot,\Delta,
      \varepsilon,S)$ associated to 
      $U\g$; $\cdot,\Delta,\varepsilon,S$ denote the product,
      coproduct, counit, antipode.
\item The quantum group \cite{dr2} $H_h\equiv(U_h\g,\bullet,\Delta_h,
      \varepsilon_h,S_h,\R)$.
\item An algebra isomorphism\cite{dr3} 
      $\varphi_h:U_h\g\rightarrow U\g[[h]]$,
      $\varphi_h\circ\bullet=\cdot\circ(\varphi_h\ot \varphi_h)$.
\item A corresponding Drinfel'd twist\cite{dr3} 
      $\F\equiv\F^{(1)}\!\ot\!\F^{(2)}\!=\!\1^{\ot^2}\!\!
      +\!O(h)\in U\g\![[h]]^{\ot^2}$:
      \[
      (\varepsilon\ot \id)\F=\1=(\id\ot \varepsilon)\F,
      \qquad\: \:\Delta_h(a)=(\varphi_h^{-1}\ot \varphi_h^{-1})\big
      \{\F\Delta[\varphi_h(a)]\F^{-1}\big\}.
      \]
\item $\gamma':=\F^{(2)}\cdot S\F^{(1)}$ and 
      $\gamma:=S\F^{-1(1)}\cdot \F^{-1(2)}$.
\item The generators $a^+_i,a^i$ of a ordinary Weyl 
      or Clifford algebra $\A$.
\item The action $\trc:U\g\times\A\rightarrow \A$; 
      $\A$ is a left module algebra under $\trc$.
\item The 
      representation $\rho$ of \g to which $a^+_i,a^i$ belong:
      \[
       x\trc a^+_i=\rho(x)^j_ia^+_j\qquad\qquad
      x\trc a^i=\rho(Sx)_j^ia^j.
      \]
\item The Jordan-Schwinger algebra homomorphism 
      $\sigma:U\g\in\A[[h]]$:
      \[
      \sigma(x):=
      \rho(x)^i_ja^+_ia^j\qquad\mbox{if~}x\in\g\qquad\qquad\qquad
      \sigma(yz)=\sigma(y)\sigma(z)
      \]
\item The generators $\tilde A^+_i,\tilde A^i$ of a deformed Weyl 
      or Clifford algebra $\Ah$.
\item The action $\tr:U_h\g\times\Ah\rightarrow \Ah$; $\Ah$ is a 
      left module algebra under $\tr$.
\item The representation $\rho_h=\rho\circ \varphi_h$
      of $U_h\g$ to which $\tilde A^+_i,\tilde A^i$ belong:
      \[
       X\tr \tilde A^+_i=\tilde\rho(X)^j_i\tilde A^+_j\qquad\qquad
      X\tr \tilde A^i=\tilde\rho(S_h X)_j^i\tilde A^j.
      \]
\item $*$-structures $*,*_h,\star,\star_h$ in $H,H_h,\A,\Ah$, if any.

\end{enumerate}

\section{Constructing the deformed generators}
\label{con}

\begin{prop}\cite{fiojmp}
One can realize the quantum group action $\tr$ on $\A[[h]]$ by setting
for any $X\in U_h\g$ and $\beta \in\A[[h]]$
(with $X_{(\bar 1)}\otimes X_{(\bar 2)}:=\Delta_h(X)$)
\begin{equation}
X\tr \beta := \sigma[\varphi_h(X_{(\bar 1)})]\,\beta\,
\sigma[\varphi_h(S_hX_{(\bar 2)})].
\end{equation}
\end{prop}

\begin{prop}\cite{fiojmp,fiocmp}
For any \g-invariants $u,v\in\A[[h]]$ the elements of $\A[[h]]$ 
\begin{equation}
\begin{array}{lll}
A_i^+ &:= & u\,\sigma(\F^{(1)})\,a_i^+\,
\sigma(S\F^{(2)}\gamma)u^{-1} \nonumber\\
A^i &:= &v\,\sigma(\gamma'S\F^{-1(2)})\,a^i\,
\sigma(\F^{-1(1)})v^{-1}.                   
\end{array}
\end{equation}
transform under $\tr$ as $\tilde A^+_i,\tilde A^i$.
\end{prop}
A suitable choice of $uv^{-1}$ may make $A^+_i,A^j$ fulfil also the 
QCR of $\Ah$ \cite{fiocmp}.  In particular we have shown the

\begin{prop}\cite{fiocmp}
If $\rho$ is the defining representation of \g,
$A^+_i,A^j$ fulfil the corresponding QCR provided
\begin{equation}
\begin{array}{llll}
uv^{-1} & = & \frac{\Gamma(n\!+\!1)}{\Gamma_{q^2}(n\!+\!1)}\qquad
\qquad & \mbox{\rm if ~}\g=sl(N) \cr
uv^{-1} & = & \frac{\Gamma[\frac 12(n\!+\!1\!+\!\frac N2-l)]
\Gamma[\frac 12(n\!+\!1\!+\!\frac N2\!+\!l)]}{\Gamma_{q^2}
[\frac 12(n\!+\!1\!+\!\frac N2\!+\!l)]
\Gamma_{q^2}[\frac 12(n\!+\!1\!+\!\frac N2-l)]}
\qquad\qquad &  \mbox{\rm if ~}\g=so(N),
\end{array}
\nonumber
\end{equation}
where $\Gamma,\Gamma_{q^2}$ are Euler's gamma-function and its
$q$-deformation,
$n:=a^+_ia^i$, $l:=\sqrt{\sigma({\cal C}_{so(N)})}$, and
${\cal C}_{so(N)}$ is the quadratic Casimir of $so(N)$.
\end{prop}

If $A^+_i,A^j$ fulfil the QCR, then also 
\begin{equation}
A^+_{i,\alpha}:=\alpha \,A^+_i\, \alpha^{-1}\qquad\qquad
A^{i,\alpha}:=\alpha \,A^i\, \alpha^{-1}
\end{equation}
will do, for any $\alpha\in\A[[h]]$ of the form $\alpha={\bf 1}+O(h)$.
By cohomological arguments one can prove
that there are no more elements in $\A[[h]]$ which do \cite{fiocmp}. 
$A^+_{i,\alpha},A^{i,\alpha}$ transform as 
$\tilde A^+_i,\tilde A^i$ under the following modified realization of
$\tr$:
\begin{equation}
X\tr^{\alpha} \beta := \alpha \sigma[\varphi_h(X_{(\bar 1)})]
\alpha^{-1}\,\beta\,
\alpha \sigma[\varphi_h(S_hX_{(\bar 2)})]\alpha^{-1}.
\end{equation}
The algebra homomorphism 
$f_{\alpha}:\Ah\rightarrow \A[[h]]$ such that
$f_{\alpha}(\tilde A^+_i)=A^+_{i,\alpha}$ and
$f_{\alpha}(\tilde A^i)=A^{i,\alpha}$ is what is usually
called a ``$q$-deforming map''.

For a compact section of $U\g$ one can choose a unitary \F,
$\F^{*\ot *}=\F^{-1}$. Then the $U\g$-covariant $*$-structure 
$(a^i)^{\star}=a^+_i$ in $\A$
is also $U_h\g$-covariant in $\A[[h]]$ and
has the form
$(A^{i,\alpha})^{\star}=A^+_{i,\alpha}$, provided we choose 
$u=v^{-1}$ and $\alpha^{\star}=\alpha^{-1}$. More formally,
under this assumption $\star\circ f_{\alpha}=f_{\alpha}\circ \star_h$,
with $\star_h$ defined by $(\tilde A^i)^{\star_h}=\tilde A^+_i$

If $H_h$ is instead a {\it triangular} deformation
of $U\g$, the previous construction can be equally performed
and leads essentially to the same results \cite{fiojmp},
 provided we choose
in the previous formulae $u\equiv v\equiv {\bf 1}$. This follows
from the triviality of the coassociator \cite{dr3}, that
characterizes triangular deformations $H_h$.

\section*{Acknowledgments}

It is a pleasure to thank J.\ Wess for
his stimulus, support and 
warm hospitality at his Institute.
This work was supported through a TMR fellowship 
granted by the European Commission, Dir. Gen. XII for Science,
Research and Development, under the contract ERBFMBICT960921.

\end{document}